\documentclass[11pt,a4paper]{amsart}
\usepackage{amscd,amssymb,amsopn,amsmath,amsthm,graphics,amsfonts,
enumerate,verbatim,calc}
\usepackage[dvips]{graphicx}
\input xy
\xyoption{all}

\usepackage{amssymb,amsmath}

\headsep=1cm \numberwithin{equation}{section}
\newtheorem{theorem}{Theorem}[section]
\newtheorem{lemma}[theorem]{Lemma}
\newtheorem{proposition}[theorem]{Proposition}
\newtheorem{corollary}[theorem]{Corollary}

\theoremstyle{definition}

\theoremstyle{remark}
\newtheorem{remark}[theorem]{Remark}
\newtheorem{example}[theorem]{Example}

\newtheorem{acknowledgement}{Acknowledgement}

\newcommand{\Ass}{\operatorname{Ass}}
\newcommand{\im}{\operatorname{im}}
\newcommand{\grade}{\operatorname{grade}}
\newcommand{\fgrade}{\operatorname{fgrade}}

\newcommand{\Spec}{\operatorname{Spec}}

\newcommand{\cd}{\operatorname{cd}}

\newcommand{\Ht}{\operatorname{ht}}
\newcommand{\pd}{\operatorname{pd}}

\newcommand{\V}{\operatorname{V}}

\newcommand{\Ext}{\operatorname{Ext}}
\newcommand{\Supp}{\operatorname{Supp}}
\newcommand{\Tor}{\operatorname{Tor}}
\newcommand{\Hom}{\operatorname{Hom}}

\newcommand{\Ann}{\operatorname{Ann}}

\newcommand{\depth}{\operatorname{depth}}

\newcommand{\Coass}{\operatorname{Coass}}

\newcommand{\vpl}{\operatornamewithlimits{\varprojlim}}

\newcommand{\lo}{\longrightarrow}
\newcommand{\fm}{\frak{m}}
\newcommand{\fp}{\frak{p}}
\newcommand{\fq}{\frak{q}}
\newcommand{\fa}{\frak{a}}
\newcommand{\fb}{\frak{b}}

\newcommand{\fn}{\frak{n}}

\begin{document}

\author[Asgharzadeh and Divaani-Aazar]{Mohsen Asgharzadeh and
Kamran Divaani-Aazar }

\title[Finiteness properties of... ]
{Finiteness properties of formal local cohomology modules and Cohen-Macaulayness}

\address{M. Asgharzadeh, Department of Mathematics, Shahid Beheshti
University, Tehran, Iran-and-School of Mathematics, Institute for Research
in Fundamental Sciences (IPM), P.O. Box 19395-5746, Tehran, Iran.}
\email{asgharzadeh@ipm.ir}
\address{K. Divaani-Aazar, Department of Mathematics, Az-Zahra
University, Vanak, Post Code 19834, Tehran, Iran-and-School of Mathematics,
Institute for Research in Fundamental Sciences (IPM), P.O. Box 19395-5746,
Tehran, Iran.} \email{kdivaani@ipm.ir}

\subjclass[2000]{13D45, 13C14, 13Exx.}

\keywords{Coassociated prime ideals, cofinite modules, Cohen-Macaulay modules, cohomological dimension, formal grade, formal local cohomology, generalized local cohomology, ring of invariants. \\
The second author was supported by a grant from IPM (No. 87130114).}

\begin{abstract} Let $\fa$ be an ideal of a local ring $(R,\fm)$
and $M$ a finitely generated $R$-module. We investigate the structure
of the formal local cohomology modules ${\vpl}_nH^i_{\fm}(M/\fa^n
M)$, $i\geq 0$. We prove several results concerning finiteness properties
of formal local cohomology modules which indicate that these modules
behave very similar to local cohomology modules. Among other things, we
prove that if $\dim R\leq 2$ or either $\fa$ is principal or
$\dim R/\fa\leq 1$, then $\Tor_j^R(R/\fa,{\vpl}_nH^i_{\fm}(M/\fa^n M))$
is Artinian for all $i$ and $j$.  Also, we examine  the notion
$\fgrade(\fa,M)$, the formal grade of $M$ with respect to $\fa$ (i.e.
the least integer $i$ such that ${\vpl}_nH^i_{\fm}(M/\fa^n M)
\neq 0$). As applications, we establish a criterion for Cohen-Macaulayness
of $M$, and also we provide an upper bound for cohomological dimension of $M$
with respect to $\fa$.
\end{abstract}

\maketitle

\section{Introduction}

Throughout this paper, all rings considered will be commutative and
Noetherian with identity and all modules are assumed to be left
unitary.  In the present paper, we investigate the structure of certain
formal cohomology modules. (For details on the notion of formal cohomology,
we refer the reader to the interesting survey article by Illusie \cite{I}.)
Let $(R,\fm)$ be a local ring, $\fa$ an ideal of $R$ and $M$ a finitely
generated $R$-module. Let $U=\Spec(R)\setminus \{\fm\}$ and $(\widehat{U},\mathcal{O}_{\widehat{u}})$ denote the formal
completion of $U$ along $\V(\fa)\setminus\{\fm\}$. Let
$\widehat{\mathcal{F}}$ denote the  $\mathcal{O}_{\widehat{u}}$
sheaf associated to ${\vpl}_nM/ \fa^n M$. Then Peskine and
Szpiro \cite[III, Proposition 2.2]{PS} have described the formal
cohomology modules $H^i(\widehat{U},\widehat{\mathcal{F}})$ via the
isomorphisms $H^i(\widehat{U},\widehat{\mathcal{F}})\cong {\vpl}_n
H^{i+1}_{\fm}(M/\fa^n M)$, $i\geq 1$.  For each
$i\geq 0$, Schenzel \cite{Sch} has called $\mathfrak{F}^i_{\fa}(M):=
{\vpl}_nH^i_{\fm}(M/\fa^n M)$ $i$th formal local cohomology module of $M$
with respect to $\fa$ and examined their structure extensively.

When $R$ is regular, Peskine and Szpiro \cite[III, Proposition 2.2]{PS}
have remarked that $\mathfrak{F}^i_{\fa}(R)\cong \Hom_R(H^{\dim R-i}_{\fa}(R),
E_R(R/\fm))$ for all $i\geq 0$.  They have used this duality result
for solving a conjecture of Hartshorne in prime characteristic, see
\cite[III, Theorem 5.1]{PS}. Also, Ogus \cite[Theorem 2.7]{O} has used
this duality result for solving this conjecture in the case $R$ contains
$\mathbb{Q}$. Since local cohomology modules of finitely generated modules
enjoy many nice finiteness and cofiniteness properties, it is rather natural
to expect that analogues of some of these properties hold for formal local
cohomology modules. In Sections 2 and 3, we obtain some finiteness properties
of formal local cohomology modules. There are some applications of these type
of finiteness results. For instance Corollary 3.7 below could be
considered as a sample. Also, \cite[Remark 5.2.1]{Hel3} and
\cite[Corollary 2.4.1]{Hel1} are two more samples of such
applications. Finding more noticeable applications (of finiteness
properties of formal local cohomology modules) will probably
need a significant effort. This could be subject of a new project
and is not adjust to the organization of this paper.

In Section 2, we deal with the question when formal local cohomology modules are
Artinian or finitely generated. We show that if an integer $t$ is such
that $\mathfrak{F}^j_{\fa}(M)$ is Artinian for all $j>t$, then
$\mathfrak{F}^t_{\fa}(M)/\fa \mathfrak{F}^t_{\fa}(M)$ is Artinian.
This immediately implies that the set $\Coass_R(\mathfrak{F}
^t_{\fa}(M))\cap \V(\fa)$ is finite and we provide an example to show that $\Coass_R(\mathfrak{F}^t_{\fa}(M))$ can be infinite, see Remark
2.8 iii) below. By \cite [Theorem 4.5]{Sch}, $l:=\dim M/\fa M$ is the largest
integer $i$ such that $\mathfrak{F}^i_
{\fa}(M)\neq 0$. Let $f$ be the least integer $i$ such that $\mathfrak{F}^i_
{\fa}(M)\neq 0$. Under some mild assumptions on $M$, Theorem 2.7 below says that
$\mathfrak{F}^f_{\fa}(M)$ and $\mathfrak{F}^l_{\fa}(M)$ are not Artinian. In view
of Theorem 2.6 iii) below, formal local cohomology modules are very seldom finitely generated.

In Section 3, we examine the Tor modules $\Tor_j^R(R/\fa,\mathfrak{F}^i_{\fa}(M))$, $i,j\geq 0$. In each of the cases a) $\dim R\leq 2$, b) $\fa$ is principal up to radical, and c) $\dim R/\fa\leq 1$, we show that $\Tor_j^R(R/\fa,\mathfrak{F}^i_{\fa}(M))$
is Artinian for all $i$ and $j$. (In particular, if $R$ is complete and $\fp$ is
a one dimensional prime ideal of $R$, then $\Tor_i^R(R/ \fp,\widehat{R_{\fp}}/R)$  is Artinian for all $i$.) Assume that either one of the above cases
holds or $R$ is regular of positive characteristic. Then,  we prove that the
Betti number $\beta^j(\fm,\mathfrak{F}^i_{\fa}(M))$ is finite for all $i$ and
$j$.

Peskine and Szpiro have introduced the notion of {\it formal grade} of an
$R$-module $N$ with respect to $\fa$ as the least integer $i$ such that ${\vpl}_nH^i_{\fm}(N/\fa^n N)\neq 0$ and denoted it by $\fgrade(\fa,N)$.
When $R$ is Gorenstein, Schenzel \cite[Lemma 4.8 d)]{Sch} has showed
that $\fgrade(\fb,R)+\cd_{\fb}(R)=\dim R$ for all ideals $\fb$ of $R$.
(Recall that for an $R$-module $N$, the cohomological dimension of $N$
with respect to an ideal $\fb$ of $R$, $\cd_{\fb}(N)$, is defined to
be the supremum of $i$'s such that $H^i_{\fb}(N)\neq 0$.) In Section 4,
we show that $$\fgrade(\fb,M)+\cd_{\fb}(M)=\dim M \  \ (\star)$$
for all ideals $\fb$ of $R$ if and only if $M$ is Cohen-Macaulay.
We provide some examples to show that the equality $(\star)$ does not
hold even for some very close generalizations of Cohen-Macaulay modules
and some very special choices of the ideal $\fb$.  Let $N$ be an
$R$-module such that $\fgrade(\fa,N)+\cd_{\fa}(N)=\dim N$ and $L$ be
a pure submodule of $N$. We investigate the question whether
the equality $\fgrade(\fa,L)+\cd_{\fa}(L)=\dim L$ holds too. We
establish two results in this direction, see Propositions 4.4 and 4.5
below.  Proposition 4.5 can be considered as a slight generalization
of the Hochster-Eagon result on Cohen-Macaulayness of invariant rings.
In \cite {Sch}, Schenzel has established several upper bounds for
$\fgrade(\fa,M)$. In particular, he \cite [Corollary 4.11]{Sch} showed
that $\fgrade(\fa,M)\leq \dim M-\cd_{\fa}(M)$. As our last result, we
establish a lower bound for $\fgrade(\fa,M)$ to the effect that
$\fgrade(\fa,M)\geq \depth M-\cd_{\fa}(M)$.

\section{Artinianness of formal local cohomology modules}

Let $\fa$ be an ideal of a local ring $(R,\fm)$ and $M$ an $R$-module. For each
integer $i\geq 0$, the $i$th formal local cohomology module of $M$ with respect
to $\fa$ is defined by $\mathfrak{F}^i_{\fa}(M):={\vpl}_nH^i_{\fm}(M/\fa^n M)$.
The {\it formal grade} of $M$ with respect to $\fa$ is defined to be the infimum of
$i$'s such that $\mathfrak{F}^i_{\fa}(M)\neq 0$ and it is denoted by $\fgrade(\fa,M)$.
(We use the usual convention that $\sup \emptyset=-\infty$ and
$\inf \emptyset=\infty$.) For two $R$-modules $M$ and $N$, the $i$th generalized local
cohomology module of $M$ and $N$ with respect to $\fa$ is defined by
$H^{i}_{\fa}(M,N):={\varinjlim}_n\Ext^{i}_{R}(M/\fa^{n}M, N)$, see \cite{Her}.
We denote the supremum of $i$'s such that $H^{i}_{\fa}(M,N)\neq 0$
by $\cd_{\fa}(M,N)$ and we will abbreviate $\cd_{\fa}(R,N)$ by $\cd_{\fa}(N)$.
Also, when $R$ is complete, the {\it canonical module} of $M$  is defined by $K_M:=\Hom_R(H_{\fm}^{\dim M}(M),E_R(R/\fm))$.

\begin{lemma} \begin{enumerate}
\item[i)]  Let $f:(T,\fn)\lo (R,\fm)$ be a ring homomorphism of local rings
such that $\fn R$ is $\fm$-primary (e.g. $R$ is integral over $T$). Let $\fb$ be
an ideal of $T$ and $M$ an $R$-module. Then for any integer $i\geq 0$, we have a natural
$R$-isomorphism $\mathfrak{F}^i_{\fb}(M)\cong \mathfrak{F}^i_{\fb R}(M)$.
\item[ii)] Let $\fa$ be an ideal of a Cohen-Macaulay complete local ring
$(R,\fm)$ and $M$ a finitely generated $R$-module. Let $K_R$ be the canonical
module of $R$. Then $$\mathfrak{F}^i_{\fa}(M)\cong\Hom_R(H^{\dim R-i}_{\fa}(M,K_R),
E_R(R/\fm))$$ for all $i\geq 0$. In particular,
$\fgrade(\fa,M)=\dim R-\cd_{\fa}(M, K_R)$.
\end{enumerate}
\end{lemma}

{\bf Proof.} i) This is an immediate  consequence of the Independence Theorem for local
cohomology modules, see \cite [Theorem 4.2.1]{BS}.

ii) The proof of the existence of these isomorphisms is the same as the proof of \cite[Remark 3.6]{Sch}, however for the sake of completeness
we include it here. For each integer $i\geq 0$, Grothendieck's Local Duality Theorem
\cite[Theorem 11.2.8]{BS} yields the isomorphism $H^i_{\fm}(M/\fa^n M)\cong
\Hom_R(\Ext_R^{\dim R-i}(M/\fa^n M,K_R),E_R(R/\fm))$ for all $n\geq 0$. Thus
$$\begin{array}{ll} \mathfrak{F}^i_{\fa}(M)&=\underset{n}{\vpl}H^i_{\fm}(M/\fa^n
M) \\&
\cong \underset{n}{\vpl}\Hom_R(\Ext_R^{\dim R-i}(M/\fa^n M,K_R),E_R(R/\fm)) \\ &
\cong \Hom_R(H^{\dim R-i}_{\fa}(M,K_R),E_R(R/\fm)).
\end{array}$$
Next, we have
$$\begin{array}{ll} \inf\{i:\mathfrak{F}^i_{\fa}(M)\neq 0\}&
=\inf\{i:H^{\dim R-i}_{\fa}(M,K_R)\neq 0\} \\&
=\inf\{\dim R-j:H^j_{\fa}(M,K_R)\neq 0\} \\&
=\dim R-\cd_{\fa}(M,K_R).
\end{array}
$$ This completes the proof of ii). $\Box$

To prove Theorem 2.4 below, we need a couple of lemmas.

\begin{lemma} Let $\fa$ be an ideal of a local ring $(R,\fm)$ and  $M$
a finitely generated $R$-module of dimension $d$. Then $\mathfrak{F}^d_{\fa}(M)$
is Artinian.
\end{lemma}

{\bf Proof.} Let $\fa:=(x_1,\cdots,x_n)$. We argue by induction on $n$.
Let $n=1$. By \cite[Corollary 3.16]{Sch}, one has the exact sequence
$$\cdots \lo H^d_{\fm}(M) \lo \mathfrak{F}_{\fa}^d(M)\lo
\Hom_R(R_{x_1},H^{d+1}_{\fm}(M))\lo \cdots .$$ Since, by
Grothendieck's Vanishing Theorem $H^{d+1}_{\fm}(M)=0$, it turns out
that $\mathfrak{F}_{\fa}^d(M)$ is a homomorphic image of the Artinian
$R$-module $H^d_{\fm}(M)$. Now, assume that the claim holds for $n-1$ and set
$\fb:=(x_1,\cdots,x_{n-1})$. Then \cite[Theorem 3.15]{Sch}, provides
the following long exact sequence
$$\cdots \lo \mathfrak{F}_{\fb}^d(M)\lo \mathfrak{F}_{\fa}^d(M)\lo
\Hom_R(R_{x_n},\mathfrak{F}_{\fb}^{d+1}(M))\lo \cdots .$$

By \cite[Theorem 4.5]{Sch}, one has $\mathfrak{F}_{\fb}^{d+1}(M)=0$, and so
$\mathfrak{F}_{\fa}^d(M)$ is a homomorphic image of $\mathfrak{F}_{\fb}^d(M)$.
Therefore, the induction hypothesis yields that $\mathfrak{F}_{\fa}^d(M)$
is Artinian. $\Box$

In the remainder of this section, we will use the following lemma.
Among other things, it says that if $\fa$ is an ideal of a local ring $(R,\fm)$
and $M$ a finitely generated $\fa$-torsion $R$-module (i.e  $H^0_{\fa}(M)=M$), then $\mathfrak{F}_{\fa}^i(M)$
is Artinian for all $i\geq 0$.

\begin{lemma} Let $\fa$ be an ideal of a local ring $(R,\fm)$, $M$ a finitely
generated $R$-module and $N$ a submodule of $M$ which is supported in $\V(\fa)$.
Then, there is a natural isomorphism
$\mathfrak{F}_{\fa}^i(N)\cong H^{i}_{\fm}(N)$ for all $i$, and so there
exists a long exact sequence $$\cdots \lo
H^{i}_{\fm}(N)\lo \mathfrak{F}_{\fa}^i(M) \lo \mathfrak{F}_{\fa}^i(M/N)\lo
H^{i+1}_{\fm}(N)\lo \cdots .$$
\end{lemma}

{\bf Proof.} Since $\Supp_R N \subseteq \V(\fa)$, it turns out that
$N$ is annihilated by some power of $\fa$. So
$$\mathfrak{F}_{\fa}^{i}(N)\cong \underset{n}{\vpl}H^{i}_{\fm}(N/ \fa ^n
N)\cong\underset{n}{\vpl}H^{i}_{\fm}(N)\cong H^{i}_{\fm}(N)$$ for
all $i$. Next, the existence of the mentioned exact sequence is immediate,
because by \cite[Theorem 3.11]{Sch}, the short exact sequence
$0\lo N \lo M\lo M/N\lo 0$ implies the long exact sequence $$\cdots \lo
\mathfrak{F}_{\fa}^i(N)\lo \mathfrak{F}_{\fa}^i(M) \lo \mathfrak{F}_{\fa}^i(M/N)\lo
\mathfrak{F}_{\fa}^{i+1}(N)\lo \cdots.  \ \  \Box$$

Now, we are in the position to present our first main result.

\begin{theorem} Let $\fa$ be an ideal of a local ring $(R,\fm)$ and
$M$ a finitely generated $R$-module. Assume that the integer $t$ is
such that $\mathfrak{F}^i_{\fa}(M)$ is Artinian for all $i>t$. Then
$\mathfrak{F}^t_{\fa}(M)/\fa \mathfrak{F}^t_{\fa}(M)$ is Artinian.
\end{theorem}

{\bf Proof.}  We use induction on $n:=\dim M$. For $n=0$, we have
$\mathfrak{F}_{\fa}^{i}(M)=0$ for all $i>0$ and $\mathfrak{F}_{\fa}^{0}(M)$ is
Artinian by Lemma 2.2. So, in this case the claim holds. Now, let
$n>0$ and assume that the claim holds for all values less than $n$.
By Lemma 2.3, one has the following long exact sequence
$$\cdots \lo H^i_{\fm}(\Gamma_{\fa}(M)) \lo \mathfrak{F}_{\fa}^i(M)
\lo \mathfrak{F}_{\fa}^i(M/\Gamma_{\fa}(M))\lo
H^{i+1}_{\fm}(\Gamma_{\fa}(M))\lo \cdots .$$ So
$\mathfrak{F}_{\fa}^i(M/\Gamma_{\fa}(M))$ is Artinian for all $i>t$. We
split the exact sequence $$ H^t_{\fm}(\Gamma_{\fa}(M)) \lo
\mathfrak{F}_{\fa}^t(M)\stackrel{\varphi}\lo
\mathfrak{F}_{\fa}^t(M/\Gamma_{\fa}(M))\stackrel{\psi}\lo
H^{t+1}_{\fm}(\Gamma_{\fa}(M))$$ to the exact sequences
$$0\longrightarrow \ker \varphi \longrightarrow \mathfrak{F}^{t}_{\fa}(M)
\longrightarrow\im \varphi \longrightarrow 0 $$ and $$ 0\lo \im
\varphi \longrightarrow
\mathfrak{F}^{t}_{\fa}(M/\Gamma_{\fa}(M))\longrightarrow \im
\psi\longrightarrow 0.$$ From these exact sequences, we deduce the
following exact sequences
$$\frac{\ker \varphi}{\fa \ker \varphi}\lo
\frac{\mathfrak{F}^{t}_{\fa}(M)}{\fa \mathfrak{F}^{t}_{\fa}(M)}\lo \frac{\im
\varphi} {\fa \im \varphi} \ \  (\star)$$
 and $$Tor^{R}_{1}(R/\fa,\im \psi)\lo \frac{\im \varphi}{\fa \im \varphi}\lo
\frac{\mathfrak{F}^{t}_{\fa}(M/\Gamma_{\fa}(M))}{\fa
\mathfrak{F}^{t}_{\fa}(M/\Gamma_{\fa}(M))}. \   \ (\star,\star)$$
Since $\ker \varphi$ and $\im \psi$ are Artinian, in view of $(\star)$
and $(\star,\star)$, it turns out that if
$\frac{\mathfrak{F}^{t}_{\fa}(M/\Gamma_{\fa}(M))}{\fa
\mathfrak{F}^{t}_{\fa}(M/\Gamma_{\fa}(M))}$ is Artinian, then
$\frac{\mathfrak{F}^{t}_{\fa}(M)}{\fa \mathfrak{F}^{t}_{\fa}(M)}$ is also
Artinian. So, we may and do assume that $M$ is $\fa$-torsion free. Take
$x\in\fa\setminus \bigcup _{\fp\in \Ass_RM}\fp$. Then $\dim M/xM=n-1$.
By \cite[Theorem 3.11]{Sch}, the exact sequence
$0\longrightarrow M \stackrel{x}\longrightarrow
M\longrightarrow M/xM\longrightarrow 0$ implies the following long exact
sequence of formal local cohomology modules
$$\cdots \lo \mathfrak{F}_{\fa}^i(M)\stackrel{x}\lo \mathfrak{F}_{\fa}^i(M)\lo
\mathfrak{F}_{\fa}^i(M/xM)\lo \mathfrak{F}_{\fa}^{i+1}(M)\lo \cdots .$$ It
yields that $\mathfrak{F}^{j}_{\fa}(M/xM)$ is  Artinian for all $j>t$.
Thus $\frac{\mathfrak{F}^{t}_{\fa}(M/xM)}{\fa \mathfrak{F}^{t}_{\fa}(M/xM)}$ is
Artinian by the induction hypothesis. Now, consider the exact sequence $$
\mathfrak{F}^{t}_{\fa}(M)\stackrel{x}\longrightarrow
\mathfrak{F}^{t}_{\fa}(M)\stackrel{f}\longrightarrow
\mathfrak{F}^{t}_{\fa}(M/xM)\stackrel{g}\longrightarrow
\mathfrak{F}^{t+1}_{\fa}(M),$$ which induces the exact
sequences $$0\longrightarrow \im f \longrightarrow
\mathfrak{F}^{t}_{\fa}(M/xM) \longrightarrow\im g \longrightarrow 0$$ and
$$\mathfrak{F}^{t}_{\fa}(M)\stackrel{x}\longrightarrow
\mathfrak{F}^{t}_{\fa}(M)\longrightarrow \im f\longrightarrow 0.$$
Therefore, we can obtain the following two exact sequences
$$Tor^{R}_{1}(R/\fa,\im g )\lo \frac{\im f}{\fa \im f}\lo
\frac{\mathfrak{F}^{t}_{\fa}(M/xM)}{\fa \mathfrak{F}^{t}_{\fa}(M/xM)},$$ and
$$\frac{\mathfrak{F}^{t}_{\fa}(M)}{\fa \mathfrak{F}^{t}_{\fa}(M)}\stackrel{x}
\longrightarrow \frac{\mathfrak{F}^{t}_{\fa}(M)}{\fa
\mathfrak{F}^{t}_{\fa}(M)}\longrightarrow \frac{\im f}{\fa \im
f}\longrightarrow0.$$ Since $x\in \fa$, from the later exact sequence,
we get that $\frac{\im f}{\fa \im f}\cong
\frac{\mathfrak{F}^{t}_{\fa}(M)}{\fa \mathfrak{F}^{t}_{\fa}(M)}$. Now, since
$Tor^{R}_{1}(R/\fa,\im g)$ and $\frac{\mathfrak{F}^t_{\fa}(M/xM)}{\fa
\mathfrak{F}^t_{\fa}(M/xM)}$ are Artinian, the claim follows. $\Box$

The statement of the corollary below involves the notion of coassociated
prime ideals. For convenient of the reader, we review this notion
briefly in below. For an $R$-module $X$, a prime ideal $\fp$ of $R$ is said
to be a {\it coassociated prime ideal} of $X$ if there exists
an Artinian quotient $Y$ of $X$ such that $\fp=\Ann_RY$. The set of all coassociated
prime ideals of $X$ is denoted by $\Coass_RX$. It is clear from the definition
that the set of coassociated prime ideals of any quotient of $X$ is contained
in $\Coass_RX$. It is known that if $X$ is Artinian, then the set $\Coass_RX$
is finite.  For more details on the notion of coassociated prime ideals,
we refer the reader to e.g. \cite{DT}. Let $\fa$ be an ideal of a ring $R$.
One can easily check that for any $R$-module $X$, $\Coass_RX/\fa X=\Coass_RX\cap \V(\fa)$.
So, we record the following immediate corollary.

\begin{corollary}  Let $\fa$ be an ideal of a local ring $(R,\fm)$ and
$M$ a finitely generated $R$-module. Assume that the integer $t$ is
such that $\mathfrak{F}^i_{\fa}(M)$ is Artinian for all $i>t$. Then
$\Coass_R(\mathfrak{F}^{t}_{\fa}(M))\cap \V(\fa)$ is finite.
\end{corollary}

The next result indicates that formal local cohomology modules are very
seldom finitely generated.

\begin{theorem} Let $\fa$ be an ideal of a local ring $(R,\fm)$ and
$M$ a finitely generated $R$-module. Then the
following assertions hold.
\begin{enumerate}
\item[i)] $\mathfrak{F}^{0}_{\fa}(M)$ is a finitely generated
$\widehat{R}$-module. In addition, if $\dim M/\fa M=0$, then $\mathfrak{F}^{0}_{\fa}(M)\cong \widehat{M}$.
\item[ii)] Assume that $\ell:=\dim M/\fa M>0$. Then $\mathfrak{F}^{\ell}_{\fa}(M)$
is not a finitely generated $R$-module.
\item[iii)] Assume that $R$ is a regular ring containing a field.
Then for any integer $i$,  the $R$-module $\mathfrak{F}^{i}_{\fa}(R)$ is
free whenever it is finitely generated.
\item[iv)] Let $t<\depth_RM$ be an integer such that
$\mathfrak{F}^{i}_{\fa}(M)$ is finitely generated for all $i<t$. Then
$\Hom_R(R/\fa,\mathfrak{F}^{t}_{\fa}(M))$ is a finitely generated
$\widehat{R}$-module.
\end{enumerate}
\end{theorem}

{\bf Proof.} i) Since $\mathfrak{F}^{0}_{\fa}(M)\cong\mathfrak{F}^{0}_{\fa
\widehat{R}}(\widehat{M})$ and $\dim_R(M/ \fa
M)=\dim_{\widehat{R}}(\widehat{M}/(\fa \widehat{R})\widehat{M})$, we
may assume that $M$ is complete in $\fm$-adic topology. So, $M$
is also complete in $\fa$-adic topology.  Hence $$\mathfrak{F}^{0}_{\fa}(M)=
\underset{n}{\vpl}H^0_{\fm}(M/\fa^n M)\subseteq \underset{n}{\vpl}(M/\fa^n M)=M.$$
Now, assume that $\dim M/\fa M=0$. Then for any integer $n\geq 0$, the module
$M/\fa^n M$ is Artinian, and so $H^0_{\fm}(M/\fa^n M)=M/\fa^n M$. Therefore $\mathfrak{F}^{0}_{\fa}(M)=\underset{n}{\vpl}(M/\fa^n M)=M$.

ii) Since $M/\fa M$ is $\fa$-torsion, Lemma 2.3 yields that
$\mathfrak{F}^{\ell}_{\fa}(M/\fa M)\cong  H_{\fm}^{\ell}(M/\fa M)$. Hence by \cite
[Theorem 3.11]{Sch}, from the short exact sequence $$0\lo \fa M \lo M\lo M/ \fa M
\lo 0,$$ we deduce the exact sequence
$\mathfrak{F}^{\ell}_{\fa}(M) \lo H_{\fm}^{\ell}(M/ \fa M) \lo
\mathfrak{F}^{\ell+1}_{\fa}(\fa M)$. Since $\fa M/\fa^2M$ is annihilated by
$\fa$, one concludes that $\fa M/\fa^2 M$ is supported in $\Supp_RM \cap
\V(\fa)=\Supp(M/\fa M),$ and so $\dim (\fa M/\fa^2M)\leq \dim(M/\fa M)$.
This yields $\mathfrak{F}^{\ell+1}_{\fa}(\fa M)=0$, by \cite [Theorem 4.5]{Sch}.
Therefore, since by \cite [Remark 2.5]{Hel2}, the $R$-module
$H_{\fm}^{\ell}(M/\fa M)$ is not finitely generated,
$\mathfrak{F}^{\ell}_{\fa}(M)$ cannot be finitely generated.

iii) Let $i\geq 0$ be an integer such that $\mathfrak{F}^i_{\fa}(R)$ is
a finitely generated $R$-module. Then by Lemma 2.1 ii), one has $$\mathfrak{F}
^i_{\fa}(R)\cong \mathfrak{F}^i_{\fa \widehat{R}}(\widehat{R})\cong \Hom_{\widehat{R}}(H^{\dim R-i}_{\fa \widehat{R}}(\widehat{R}),
E_R(R/\fm)).$$ So, the Matlis duality implies that $H^{\dim R-i}_{\fa \widehat{R}}(\widehat{R})$ is an Artinian $\widehat{R}$-module. Hence
$H^{\dim R-i}_{\fa \widehat{R}}(\widehat{R})$ is an injective
$\widehat{R}$-module, see \cite [Corollary 3.8]{HS} for the positive
characteristic case and \cite [Corollary 3.6 b)]{L} for the other case. Thus,
$H^{\dim R-i}_{\fa \widehat{R}}(\widehat{R})$ is a direct sum of finitely
many copies of $E_R(R/\fm)$, and so $\mathfrak{F}^i_{\fa}(R)$ is a finitely
generated flat $R$-module. This yields the conclusion, because any
finitely generated flat $R$-module is free.

iv) One has $\mathfrak{F}^{t}_{\fa}(M)\cong \mathfrak{F}^{t}_{\fa \widehat{R}}
(\widehat{M})$, and so
$$\begin{array}{ll}
\Hom_{\widehat{R}}(\widehat{R}/\fa\widehat{R}, \mathfrak{F}^{t}_{\fa \widehat{R}}
(\widehat{M}))
&\cong\Hom_{\widehat{R}}(R/\fa\otimes_R\widehat{R},\mathfrak{F}^{t}_{\fa}(M))\\
&\cong\Hom_R(R/\fa,\Hom_{\widehat{R}}(\widehat{R},\mathfrak{F}^{t}_{\fa}(M)))\\
&\cong\Hom_R(R/\fa,\mathfrak{F}^{t}_{\fa}(M)).\\
\end{array}
$$
Hence, we  can assume that $R$ is complete. By Cohen's Structure Theorem,
there exists a complete regular local ring $(T,\fn)$ such
that $R\cong T/J$ for some ideal $J$ of $T$. Set $b:=\fa\cap T$. Then by
Lemma 2.1 i), $\mathfrak{F}^i_{\fa}(M)\cong \mathfrak{F}^i_{\fb}(M)$ for all $i\geq 0$.
Also, the two $R$-modules $\Hom_R(R/\fa,\mathfrak{F}^{t}_{\fa}(M))$ and $\Hom_T(T/\fb,\mathfrak{F}^{t}_{\fb}(M))$ are isomorphic and
$\depth_TM=\depth_RM$. For any $R$-module $X$, being finitely generated
as an $R$-module is the same as being finitely generated as a $T$-module.
Thus we may and do assume that $R=T$. Let $d:=\dim R$. Then
by Lemma 2.1 ii), $\mathfrak{F}^{i}_{\fa}(M)\cong\Hom_R(H^{d-i}_{\fa}(M,R),E_R(R/\fm))$
for all $i$, and so it follows that $H^{j}_{\fa}(M,R)$ is Artinian for all
$j>d-t$. On the other hand, by the Auslander-Buchsbaum
formula, $\pd_RM=\dim R-\depth_RM<d-t$. So
\cite [Theorem 3.1]{ADT} yields that $H^{d-t}_{\fa}(M,R)/\fa H^{d-t}_{\fa}(M,R)$
is Artinian. Thus $$\Hom_R(R/\fa,\mathfrak{F}^{t}_{\fa}(M))\cong\Hom_R(H^{d-t}_
{\fa}(M,R)/\fa H^{d-t}_{\fa}(M,R),E_R(R/\fm))$$ is finitely generated,
as required. $\Box$

Part i) of the following result asserts that in Theorem 2.6 iv) if $t=\fgrade(\fa,M)$,
then $\mathfrak{F}^t_{\fa}(M)$ is not Artinian.

\begin{theorem} Let $\fa$ be an ideal of a local ring $(R,\fm)$.
\begin{enumerate}
\item[i)] If $M$ is a finitely generated $R$-module such that
$f:=\fgrade(\fa,M)<\depth M$, then $\mathfrak{F}^f_{\fa}(M)$
is not Artinian.
\item[ii)] If $R$ is Cohen-Macaulay and $\Ht\fa>0$, then
$\mathfrak{F}_{\fa}^{\dim R/\fa}(R)$ is not Artinian.
\end{enumerate}
\end{theorem}

{\bf Proof.} Without loss of generality we can assume that $R$ is complete.

i) By using Cohen's Structure Theorem, there exists a complete regular local ring
$(T,\fn)$ such that $R\cong T/J$ for some ideal $J$ of $T$. Set $b:=\fa\cap T$.
Then by Lemma 2.1, one has $$\mathfrak{F}^f_{\fa}(M)\cong
\mathfrak{F}^f_{\fb}(M)\cong\Hom_T(H^{\dim T-f}_{\fb}(M,T),E_T(T/\fn))$$ and
$(c:=)\dim T-f=\cd_{\fb}(M,T)$.  By induction on $\pd_TL$, it is easy to see
that $H^i_{\fb}(M,L)=0$ for all $T$-modules $L$ and all $i>c$. (Note that if $0\lo X \lo Y\lo Z\lo 0$ is an exact sequence of $T$-modules and $T$-homomorphisms, then one has the exact sequence $$\cdots \lo H^i_{\fb}(M,Y)\lo H^i_{\fb}(M,Z)\lo H^{i+1}_{\fb}(M,X)\lo \cdots .)$$ Thus the functor $H^{c}_{\fb}(M,-)$ is right exact. So
$$
\begin{array}{ll} \frac{H^{c}_{\fb}(M,T)}{\fb H^{c}_{\fb}(M,T)}&\cong
H^{c}_{\fb}(M,T/ \fb)
\\&\cong \Ext_T^c(M,T/\fb).
\\
\end{array}
$$
Note that since $T/\fb$ is $\fb$-torsion, \cite [Corollary 2.8 i)]{DH} implies
that $H^{c}_{\fb}(M,T/ \fb)\cong \Ext_T^c(M,T/\fb)$. Now, by the Auslander-Buchsbaum
formula, we have $$\pd_TM=\dim T-\depth_TM<\dim T-f=c,$$ and so $\frac{H^{c}
_{\fb}(M,T)}{\fb H^{c}_{\fb}(M,T)}=0$. If $H^{c}_{\fb}(M,T)$ is
finitely generated, then Nakayama's Lemma implies that $H^{c}_{\fb}(M,T)=0$,
which is a contradiction. Therefore $H^{c}_{\fb}(M,T)$ is not a finitely generated
$T$-module. Hence $\mathfrak{F}^f_{\fb}(M)$ is not an Artinian $T$-module, and so
$\mathfrak{F}^f_{\fa}(M)$ is not an Artinian $R$-module.

ii) By Lemma 2.1 ii), one has $\mathfrak{F}^i_{\fa}(R)\cong\Hom_R(H^{\dim R-i}_{\fa}(K_R),E_R(R/\fm))$ for all
$i\geq 0$. Hence the supremum of the integers $i$ such that $\mathfrak{F}^i_{\fa}(R)$
is not Artinian is equal to $\dim R-f_{\fa}(K_R)$. (Recall that $f_{\fa}(K_R)$ denotes the infimum of the integers $i$ such that $H^i_{\fa}(K_R)$ is not finitely generated.)
By \cite [Remark 2.7 ii)]{ADT}, we know that $f_{\fa}(K_R)=\Ht_{K_R}(\fa)=\Ht\fa$. Thus $\dim R/\fa$ is
the supremum of the integers $i$ such that $\mathfrak{F}^i_{\fa}(R)$ is not Artinian.
$\Box$

The following remark indicates that both the assumptions and the assertions of our
results in this section are sharp.

\begin{remark} Let $\fa$ be an ideal of a local ring $(R,\fm)$ and $M$ a
finitely generated $R$-module.
\begin{enumerate}
\item[i)] By \cite [Lemma 4.8 b)]{Sch}, one has $\fgrade(\fa,M)\leq \depth M$.
So, the condition $\fgrade(\fa,M)<\depth M$ is not a big assumption in Theorem 2.7 i). However, it cannot be dropped. To realize this, assume that $M$ is a nonzero $R$-module of finite length. Then $\mathfrak{F}^0_{\fm}(M)\cong M$ is an Artinian $R$-module.  Note that $\fgrade(\fm,M)=\depth M=0$.
\item[ii)] By Lemma 2.2, $\mathfrak{F}_{\fa}^{\dim M}(M)$ is Artinian.
Schenzel \cite[Theorem 4.5]{Sch} has showed that $\ell:=\dim M/\fa M$ is the
largest integer $i$ such that $\mathfrak{F}_{\fa}^i(M)\neq 0$. It is natural to ask
whether $\mathfrak{F}_{\fa}^{\ell}(M)$ is Artinian. However, Theorem 2.7 ii) shows that in the case $R$ is Cohen-Macaulay, the module $\mathfrak{F}_{\fa}^{\dim R/\fa}(R)$
is Artinian if and only if $\Ht \fa=0$. So, easily one can construct an example such that $\mathfrak{F}_{\fa}^{\ell}(M)$ is not Artinian.
\item[iii)] By Corollary 2.5, we know that if for an integer $t$, all the formal
local cohomology modules $\mathfrak{F}_{\fa}^{t+1}(M), \mathfrak{F}_{\fa}^{t+2}(M),
\dots $ are Artinian, then $\Coass_R(\mathfrak{F}_{\fa}^t(M))\cap \V(\fa)$ is finite.
It is rather natural to ask whether $\Coass_R(\mathfrak{F}_{\fa}^t(M))$ is also finite. This is not the case. To this end, let $T:=\mathbb{Q}[X,Y]_{(X,Y)}$ and $\fa:=(X,Y)T$. Then $\mathfrak{F}^0_{\fa}(T)=\widehat{T}=\mathbb{Q}[[X,Y]]$ and $\mathfrak{F}^i_{\fa}(T)=0$ for all $i>0$. For each integer $n$, let $\fp_n:=(X-nY)T$. Then it is easy to see that $T/\fp_n\cong \mathbb{Q}[Y]_{(Y)}$, and so it is not a complete local ring. By \cite [Beispiel 2.4]{Z}, $\Coass_T\widehat{T}=\{\fa\}
\cup \{\fp\in \Spec T:T/\fp \text { is not complete }\}.$ Hence $\Coass_T(\mathfrak{F}^0_{\fa}(T))$ is not finite.
\item[iv)] Formal local cohomology modules are pure injective. To realize this,
let $D^{\bullet}_{\widehat{R}}$ be a normalized dualizing complex of $\widehat{R}$. Then by \cite[Theorem 3.5]{Sch}, for any $i$, one has $$\mathfrak{F}^i_{\fa}(M)\cong \mathfrak{F}^i_{\fa \widehat{R}}(\widehat{M})\cong
\Hom_{\widehat{R}}(H^{-i}_{\fa \widehat{R}}(\Hom_{\widehat{R}}(\widehat{M},D^{\bullet}_{\widehat{R}})),E_R(R/ \fm)).$$ Hence $\mathfrak{F}^i_{\fa}(M)$ is a pure injective $R$-module, see Lemma 4.1 (and its preceding paragraph) in \cite {M2}. In particular, $\mathfrak{F}^i_{\fa}(M)$'s are cotorsion (i.e. $\Ext_R^j(F,\mathfrak{F}^i_{\fa}(M))=0$ for all $j\geq 1$ and all flat $R$-modules $F$).
\item[v)] One can also prove Theorem 2.4 by an argument similar to the proof
of Theorem 2.6 iv). But, we prefer the more direct existing argument.
\item[vi)] In Theorem 2.6, we have seen that the formal local cohomology
modules $\mathfrak{F}^i_{\fa}(M)$ are very seldom finitely generated.
In fact, even their set of associated primes might be infinite.
For instance, let $R$ be complete Gorenstein and equicharacteristic with
$\dim R>2$. Let $\fp$ be a prime ideal of $R$ of height $2$ and take
$x\in \fm-\fp$. Then by \cite[Corollary 2.2.2]{Hel1}, $\Ass_R(\mathfrak{F}^{\dim R-1}_{(x)}(R))=\Spec R\setminus \V((x))$. Since $\Ht \fp=2$, there are
infinitely many prime ideals of $R$ which are contained in $\fp$,
and so $\Ass_R(\mathfrak{F}^{\dim R-1}_{(x)}(R))$ is infinite.
\item[vii)] Lemma 2.1 i) can be considered as the analogue of the Independence Theorem
(for local cohomology modules) for formal local cohomology modules. One might
also expect that the analogue of the Flat Base Change Theorem (for local cohomology
modules) holds for formal local cohomology modules. Let us be more precise. Let
$f:(R,\fm)\lo (U,\fn)$ be a flat local homomorphism, $M$ a finitely generated
$R$-module and $\fa$ an ideal of $R$. Are the two $U$-modules $\mathfrak{F}^i_{\fa}(M)\otimes_RU$
and $\mathfrak{F}^i_{\fa U}(M\otimes_RU)$ isomorphic for all $i\geq 0$? This is not the
case. For example, let $k$ be a field and in the ring $R:=k[[W,X,Y,Z]]$, set $\fp:=(W,X)$
and $\fa:=\fp\cap (Y,Z)$. Then $\mathfrak{F}^1_{\fa}(R)_{\fp}\cong R_{\fp}$, see \cite
[Example 5.2]{Sch}. On the other hand $\mathfrak{F}^1_{\fa R_{\fp}}(R_{\fp})\cong
\mathfrak{F}^1_{\fp R_{\fp}}(R_{\fp})=0$.
\end{enumerate}
\end{remark}

\section{Artinianess of the modules $\Tor_j^R(R/\fa,\mathfrak{F}^i_{\fa}(M))$}

Let $\fa$ be an ideal of $R$ and $X$ an $R$-module. The module $X$ is said to be $\fa$-cofinite if it is supported in $\V(\fa)$, and $\Ext^i_R(R/ \fa,X)$
is finitely generated for all $i$. Let $M$ be a finitely generated $R$-module.
It is known that if either $\fa$ is principal or $R$ is local and $\dim R/\fa=1$,
then  the modules $H^i_{\fa}(M)$ are $\fa$-cofinite, see \cite [Theorem 1]{K} for
the case $\fa$ is principal and \cite [Theorem 1]{DM} and \cite {Y} for the other
case.  As the main results of this section, we prove
that if $\dim R\leq 2$ or either $\fa$ is principal or $\dim R/\fa\leq 1$, then
$\Tor_j^R(R/ \fa,\mathfrak{F}^i_{\fa}(M))$ is Artinian for all $i$ and $j$.

\begin{lemma} Let $\fa$ be an ideal of $R$, $X$ an $R$-module and $n\geq 0$ an integer.
Then $\Tor_i^R(R/\fa,X)$ is Artinian for all $i<n$ if and only if $\Tor_i^R(M,X)$ is
Artinian for any finitely generated $R$-module $M$ which is supported in
$\V(\fa)$ and all $i<n$.
\end{lemma}

{\bf Proof.} Using Gruson's Theorem \cite [Theorem 4.1]{V}, the proof
is an straightforward adaption of the argument of \cite [Proposition 1]{DM}. $\Box$

\begin{theorem} Let $\fa$ be an ideal of a local ring $(R,\fm)$ and $M$ a finitely
generated $R$-module. Assume that $\fa$ is principal up to radical. Then $\Tor_j^R(R/\fa,
\mathfrak{F}^i_{\fa}(M))$ is Artinian for all $i$ and $j$.
\end{theorem}

{\bf Proof.} Let $\fa$ and $\fb$ be two ideals of $R$ with the same radical
and $X$ and $Y$
two $R$-modules. Then for each integer $i$, one has $H^i_{\fa}(X,Y)\cong H^i_{\fb}(X,Y)$.
Now, the argument given in the beginning of the proof Theorem 2.7 indicates that
$\mathfrak{F}^i_{\fa}(M)\cong \mathfrak{F}^i_{\fb}(M)$. Thus in view of Lemma 3.1, without loss of generality, we may assume that
$\fa$ is principal. So, let $\fa=(x)$ for some $x\in R$ and let $i\geq 0$ be
an integer. Then by \cite[Corollary 3.16]{Sch},
there exists the following long exact sequence $$\cdots  \lo
H^i_{\fm}(M)\stackrel{f} \lo \mathfrak{F}_{\fa}^i(M)\stackrel{g}\lo
\Hom_R(R_{x},H^{i+1}_{\fm}(M))\stackrel{h}\lo H^{i+1}_{\fm}(M)\lo
\cdots .$$ Consider the following two short exact sequences

$$0\longrightarrow \im f \longrightarrow \mathfrak{F}^{i}_{\fa}(M)
\longrightarrow\im g \longrightarrow 0  \   \ (\star)$$ and $$ 0 \lo
\im g\longrightarrow \Hom_R(R_x,H^{i+1}_{\fm}(M))\longrightarrow \im
h\longrightarrow 0. \  \  (\star,\star)$$ Since $\im f$ and $\im h$
are Artinian, it turns out that $\Tor_j^R(R/\fa, \im f)$ and
$\Tor_j^R(R/\fa, \im h)$ are Artinian for all $j\geq 0$. Since the
map induced by multiplication by $x$ on $R_x$ is an isomorphism and
$x\in \fa$, we conclude that $\Tor_j^R(R/\fa,\Hom_R(R_x,H^{i+1}_{\fm}(M)))=0$
for all $j$. Thus from the long exact sequence of Tor modules which is
induced by $(\star,\star)$, it turns out that $\Tor_j^R(R/\fa, \im g)$ is Artinian
for all $j\geq 0$. Now, the long exact sequence of Tor modules which is induced
by $(\star)$ completes the proof. $\Box$

Theorem 3.6 below is our next main result. To prove it, we need the following
three lemmas. The first two lemmas enable us to reduce to the case when $R$ is a
complete regular local ring. Our approach for this task is motivated by that of
Delfino and Marley for proving their main result in \cite{DM}.

\begin{lemma} Let $f:T \lo R$ be a module-finite ring homomorphism and $X$ an
$R$-module. Then $X$ is Artinian as an $R$-module if and only if it
is Artinian as a $T$-module.
\end{lemma}

{\bf Proof.} Clearly if $X$ is Artinian as a $T$-module, then it is
also Artinian as an $R$-module. Now, assume that $X$ is Artinian as an
$R$-module. Then there are finitely many maximal ideals
$\fm_1,\cdots,\fm_t$ of $R$ such that $X$ is isomorphic to an
$R$-submodule of $\bigoplus_{i=1}^t E_R(R/\fm_i)$. So, it is enough to
prove the claim only for Artinian modules of the form $X=E_R(R/\fm)$,
where $\fm$ is a maximal ideal of $R$.
Since $f:T \lo R$ is  module-finite, it follows that $n:=\fm\cap T$
is a maximal ideal of $T$. The homomorphism $f$ induces a natural
$T$-monomorphism $f^*:T/\fn\lo R/\fm$, which in turn induces a
surjective $T$-homomorphism
$\Hom_T(R/\fm,E_T(T/\fn))\lo \Hom_T(T/\fn,E_T(T/\fn))$. This yields that
$\Hom_T(R/\fm,E_T(T/\fn))$ is nonzero. Let $Y:=\Hom_T(R,E_T(T/\fn))$. Then, it
is easy to see that $Y$ is an injective $R$-module and an Artinian $T$-module.
We have $$\Hom_R(R/\fm,Y)\cong
\Hom_T(R/\fm,E_T(T/\fn))\neq 0.$$
So $\fm\in \Ass_RY$, and hence $E_R(R/\fm)$ is a direct summand of the
injective $R$-module $Y$. Therefore, $E_R(R/\fm)$ is Artinian as a $T$-module,
as required. $\Box$

In Theorem 3.6 below, we use a special case of the following result, in which $T$ is
local and $R$ is a homeomorphic image of $T$. But here we prefer to include the following general setting for other possible applications in future.

\begin{lemma} Let $f:T\lo R$  be a module-finite ring homomorphism.
Let $\fb$ be an ideal of $T$ and $X$ an $R$-module. Then the $R$-module
$\Tor_i^R(R/\fb R,X)$ is Artinian for all $i$ if and only if the $T$-module $\Tor_i^T(T/\fb,X)$ is Artinian for all $i$.
\end{lemma}

{\bf Proof.} By \cite[Theorem 11.62]{R}, we have the following spectral
sequence $$E^2_{\fp,\fq}:=\Tor_p^R(\Tor_q^T(T/\fb,R),X)\underset{p}
\Longrightarrow \Tor_{p+q}^T(T/\fb,X).$$ First suppose that the
$R$-module $\Tor_i^R(R/\fb R,X)$ is Artinian for all $i$. For any
$q\geq 0$, the $R$-module $\Tor_q^T(T/\fb,R)$ is finitely generated
and is supported in $\V(\fb R)$. Hence Lemma 3.1, implies that
the $R$-module $E^2_{\fp,\fq}$ is Artinian for all $p,q$. For
each $n$, there exists a filtration
$$0=H_{-1}\subseteq H_0\subseteq\cdots\subseteq H_n=\Tor_{n }^T(T/\fb,X)$$
of submodules of $\Tor_{n }^T(T/\fb,X)$ such that $H_i/H_{i-1}\cong
E^{\infty}_{i,n-i}$ for all $i=0,\cdots,n$. But for each $i$, $E^{\infty}_{i,n-i}$
is a subquotient of $E^{2}_{i,n-i}$, and so it is an Artinian $R$-module.
Thus $\Tor_{n }^T(T/\fb,X)$ is an Artinian $R$-module for all $n\geq 0$.
Hence, by Lemma 3.3, $\Tor_{n }^T(T/\fb,X)$ is an Artinian $T$-module for all
$n\geq 0$.

Conversely, assume that the $T$-module $\Tor_i^T(T/\fb,X)$ is Artinian for
all $i$. By induction on $n$, we prove that $E^{2}_{n,0}\cong
\Tor_n^R(R/\fb R,X)$ is an Artinian $R$-module for all $n$. For $n=0$, one has
$E^{2}_{0,0}\cong T/\fb\otimes_TX$, so it is Artinian as a $T$-module
as well as an $R$-module. Now, assume that the claim is true for all $p<n$.
Then Lemma 3.1, implies that $E^{2}_{p,q}$ is an Artinian $R$-module
for all $p<n$ and $q\geq 0$. One has the exact sequence
$$0 \lo E^{r+1}_{n,0}\lo E^{r}_{n,0}\stackrel{d^r_{n,0}}\lo
E^{r}_{n-r,r-1}\ \ (\star)$$ for all $r\geq 2$. Since $\Tor_{n
}^T(T/\fb,X)$ is an Artinian $T$-module, it follows that
$E^{\infty}_{n,0}$ is an Artinian $T$-module. We have
$E^{\infty}_{n,0}\cong E^{r}_{n,0}$ for all $r\gg 0$. By using $(\star)$
recursively, it becomes clear that $E^{2}_{n,0}$ is an Artinian $R$-module.
$\Box$

\begin{lemma} Let $\fa$ be an ideal of a  regular complete local ring $(R,\fm)$ and
$M$ a finitely generated $R$-module.  Assume that $\dim R/\fa =1$. Then
$H^n_{\fa}(M,R)$ is $\fa$-cofinite for all $n$.
\end{lemma}

{\bf Proof.} First of all note that in view of \cite [Proposition 1]{DM}, we can and do assume that $\fa$ is radical. Hence, the assumption $\dim R/\fa=1$ yields that $\fa$
is the intersection of finitely many one dimensional prime ideals of $R$. Thus $\fa$ is unmixed, because $R$ is catenary. Let $d:=\dim R$. Since $R$ is a complete domain, the
Hartshorne-Lichtenbaum Vanishing Theorem yields that $H^d_{\fa}(R)=0$.
Thus $\grade(\fa,R)=\cd_{\fa}(R)=d-1$. So, the
spectral sequence
$$\Ext_R^p(M,H^q_{\fa}(R))\underset{p}\Longrightarrow
H_{\fa}^{q+p}(M,R)$$ collapses at $q=d-1$. Hence
$H^n_{\fa}(M,R)\cong\Ext_R^{n-d+1}(M,H^{d-1}_{\fa}(R))$ for all
$n\geq 0$. By \cite[Proposition 5.2]{B}, we know that
$$\inf\{i:H^i_{\fa}(M,R)\neq0\}= \grade(\Ann_R(M/\fa M),R)\geq
\grade(\fa,R)=d-1.$$ On the other hand, since the injective dimension of $R$ is
equal to $d$, one has $H_{\fa}^i(M,R)=0$ for all $i>d$. Hence $H^n_{\fa}(M,R)=0$ for all
$n\notin \{d,d-1\}$. By \cite[Lemma 4.7]{HK}, the
$R$-module $H^d_{\fa}(M,R)\cong \Ext_R^{1}(M,H^{d-1}_{\fa}(R))$ is
$\fa$-cofinite. Therefore, it remains to prove that
$H_{\fa}^{d-1}(M,R)\cong \Hom_R(M,H^{d-1}_{\fa}(R))$ is $\fa$-cofinite.
By \cite[Lemma 4.3]{HK}, this holds if $M$ is a submodule of a
finitely generated free $R$-module. We can construct an exact
sequence $0\lo N\lo F\lo M\lo 0,$ where $F$ is a finitely
generated free $R$-module. This short exact sequence induces the
following exact sequence
$$
\begin{array}{ll}0\lo\Hom_R(M,H^{d-1}_{\fa}(R))&\lo
\Hom_R(F,H^{d-1}_{\fa}(R))\\&
\stackrel{f}\lo\Hom_R(N,H^{d-1}_{\fa}(R))\lo \Ext^1_R(M,H^{d-1}
_{\fa}(R))\lo 0.\\
\end{array}
$$
We split it into the short exact sequences
$$0\lo\Hom_R(M,H^{d-1}_{\fa}(R))\lo\Hom_R(F,H^{d-1}_{\fa}(R)) \lo\im f\lo
0 \   \  (\star)$$ and
$$ 0\lo \im f \lo\Hom_R(N,H^{d-1}_{\fa}(R))\lo\Ext^1_R(M,H^{d-1}_{\fa}(R))
\lo0.\  \  (\star\star)$$

By \cite[Lemma 4.3]{HK}, the modules  $\Hom_R(N,H^{d-1}_{\fa}(R))$
and $\Hom_R(F,H^{d-1}_{\fa}(R))$ are $\fa$-cofinite. Since
$\Ext^1_R(M,H^{d-1}_{\fa}(R))$ and $\Hom_R(N,H^{d-1}_{\fa}(R))$ are $\fa$-cofinite,
using the long exact sequence of Ext modules that induced by $(\star\star)$, one sees
that $\im f$ is also $\fa$-cofinite. Now from $(\star)$, one concludes that $\Hom_R(M,H^{d-1}_{\fa}(R))$ is $\fa$-cofinite, as desired. $\Box$

Now, we are ready to prove the next main result of the paper.

\begin{theorem} Let $\fa$ be an ideal of a local ring $(R,\fm)$  and $M$ a
finitely generated $R$-module.  If $\dim R/\fa \leq 1$, then
$\Tor_j^R(R/\fa, \mathfrak{F}^i_{\fa}(M))$ is Artinian for all $i$ and
$j$.
\end{theorem}

{\bf Proof.} The case $\dim R/\fa=0$ is trivial. So, in below, we assume
that $\dim R/\fa=1$. Let $F_{\bullet}$ be a free resolution of the
$R$-module $R/\fa$. Then, clearly $F_{\bullet}\otimes_R\widehat{R}$ is
a free resolution of the $\widehat{R}$-module $\widehat{R}/\fa \widehat{R}$.
Hence, for any $\widehat{R}$-module $X$ and any $i\geq 0$, one has
$$\Tor_i^R(R/\fa,X)\cong H_i(F_{\bullet}\otimes_RX)\cong H_i((F_{\bullet}\otimes_R\widehat{R})\otimes_{\widehat{R}}X)\cong
\Tor_i^{\widehat{R}}(\widehat{R}/\fa \widehat{R},X).$$
Thus for any $i$ and $j$, the two $\widehat{R}$-modules $\Tor_{j}^R(R/ \fa,\mathfrak{F}^i_{\fa}(M))$ and $\Tor_{j}^{\widehat{R}}({\widehat{R}}/\fa
\widehat{R},\mathfrak{F}^i_{\fa\widehat{R}}(\widehat{M}))$ are isomorphic.
So, we may and do assume that $R$ is complete. Then by Cohen's Structure
Theorem, $R$ is a homomorphic image of a complete regular local ring
$(T,\fn)$. That is $R\cong T/J$ for some ideal $J$ of $T$. By Lemma 2.1 i) and
Lemma 3.4, we can assume that $R=T$. Since for any $R$-module $X$, one can
see easily that $\Tor_j^R(R/\fa, \Hom_R(X,E_R(R/\fm)))\cong \Hom_R(\Ext^j_R(R/\fa,X),E_R(R/\fm))$, Lemma 2.1 ii) implies that
$$\Tor_j^R(R/ \fa, \mathfrak{F}^i_{\fa}(M))\cong \Hom_R(\Ext^j_R(R/\fa,H^{d-i}_{\fa}(M,R)),E_R(R/\fm)),$$ where $d=\dim R$.
Therefore, for any $i$ and $j$, Matlis Duality asserts that
$\Tor_j^R(R/ \fa,\mathfrak{F}^i_{\fa}(M))$ is Artinian if and only if
$\Ext^j_R(R/\fa,H^{d-i}_{\fa}(M,R))$ is finitely generated. So, the claim
follows by Lemma 3.5. $\Box$

The computation that was done by Hellus in the proof of \cite [Lemma 3.2.1]{Hel3} shows
that if $\fp$ is a one dimensional prime ideal of a complete local ring $(R,\fm)$, then $\mathfrak{F}^1_{\fp}(R)\cong \widehat{R_{\fp}}/R$. Thus, one has the following corollary.

\begin{corollary} Let $\fp$ be a one dimensional prime ideal of a complete local
ring $(R,\fm)$. Then $\Tor_i^R(R/ \fp,\widehat{R_{\fp}}/R)$  is Artinian for all
$i$.
\end{corollary}

For proving the second part of our next result, we employ an argument analogue
to that used by Melkersson in \cite[Theorem 2.1]{M1}.

\begin{theorem} Let $\fa$ be an ideal of a local ring $(R,\fm)$ and
$M$ a finitely generated $R$-module. Then the following assertions
hold.
\begin{enumerate}
\item[i)] If either $i=0$ or $i=\dim M$, then $\Tor_j^R(R/\fa,\mathfrak{F}^i_{\fa}(M))$ is Artinian for all $j$.
\item[ii)] If $\dim R\leq 2$, then $\Tor_j^R(R/\fa,\mathfrak{F}^i_{\fa}(M))$
is Artinian for all $i$ and $j$.
\end{enumerate}
\end{theorem}

{\bf Proof.} i) The case $i=\dim M$ is clear by Lemma 2.2. Now, let
$i=0$. Without loss of generality, we may and do assume that $R$ is
complete. Then by Theorem 2.6 i), $\mathfrak{F}^0_{\fa}(M)$ is a finitely
generated $R$-module. \cite[Lemma 4.1]{Sch} yields that
$$\Ass_R(\mathfrak{F}^0_{\fa}(M))=\{\fp\in \Ass_RM:\dim R/(\fa+\fp)=0\},$$
and so $\V(\fa)\cap \Supp_R(\mathfrak{F}^0_{\fa}(M))\subseteq \{\fm\}$. Let $j$ be an integer. Then  $\Tor_j^R(R/\fa,\mathfrak{F}^0_{\fa}(M))$ is
supported only at $\fm$, and hence $\Tor_j^R(R/\fa,\mathfrak{F}^0_{\fa}(M))$
has finite length.

ii) As in i), we may and do assume that $R$ is complete. The cases $i=0$
and $i=2$ follow by i). Since by \cite[Theorem 4.5]{Sch},
$\mathfrak{F}^i_{\fa}(M)=0$ for all $i>2$, it remains to show that $\Tor_j^R(R/\fa,
\mathfrak{F}^1_{\fa}(M))$ is Artinian for all $j$. There are prime ideals $\fp_1,\cdots,\fp_n$ and a chain $0=M_0\subseteq M_1 \subseteq \cdots \subseteq M_n=M$
of submodules of $M$ such that $M_i/ M_{i-1}\cong R/\fp_i$ for all
$i=1,\cdots, n$. Now, we complete the argument by applying induction on $n$.
Let $n=1$, and set $A:=M\cong R/\fp_1$. By Lemma 2.1 i), we have $\mathfrak{F}^1_{\fa}(M)
\cong \mathfrak{F}^1_{\fa A}(A)$. So, in view of Lemma 3.4, it suffices
to show that $\Tor_j^A(A/\fa A,\mathfrak{F}^1_{\fa A}(A))$
is Artinian for all $j$. If $\dim A/\fa A=2$, then $\fa A=0$, and so
$\mathfrak{F}^1_{\fa A}(A)\cong H^1_{\fm A}(A)$. Therefore, the proof of
the case $n=1$ is complete. Next, assume that $n>1$ and that
the claim has been proved for $n-1$. From the short exact sequence
$$0\lo M_{n-1}\lo M\lo R/\fp_n\lo 0,$$ by \cite[Theorem 3.11]{Sch},
one has the long exact sequence
$$\begin{array}{ll}
0\rightarrow \mathfrak{F}^0_{\fa}(M_{n-1})\longrightarrow
\mathfrak{F}^0_{\fa}(M)\longrightarrow & \mathfrak{F}^0_{\fa}(R/\fp_n )
\longrightarrow  \mathfrak{F}^1_{\fa}(M_{n-1})\longrightarrow
\mathfrak{F}^1_{\fa}(M)\longrightarrow \\& \mathfrak{F}^1_{\fa}(R/\fp_n)
\longrightarrow  \mathfrak{F}^2_{\fa}(M_{n-1})\longrightarrow
\mathfrak{F}^2_{\fa}(M)\longrightarrow \mathfrak{F}^2_{\fa}(R/\fp_n)
\rightarrow 0.
\end{array}
$$
Now, by splitting this long exact sequence into short exact sequences,
we can prove that $\Tor_j^R(R/\fa,\mathfrak{F}^1_{\fa}(M))$ is
Artinian for all $j$.  Note that in a short exact sequence
$$0\lo X\lo Y\lo Z \lo 0$$ if two of the modules $X,Y$ and $Z$ have the
property that all of their $\Tor$-modules against $R/\fa$ are Artinian,
then the same property also holds for the third one. $\Box$

\begin{corollary} Let $\fa$ be an ideal of a local ring $(R,\fm)$  and $M$ a
finitely generated $R$-module.  In each of the following cases,
the Betti number $\beta^j(\fm,\mathfrak{F}^i_{\fa}(M))$ is
finite for all $i$ and $j$.
\begin{enumerate}
\item[i)] $\dim R\leq 2$.
\item[ii)] $\fa$ is principal up to radical.
\item[iii)] $\dim R/\fa\leq 1$.
\item[iv)] $R$ is regular of positive characteristic.
\end{enumerate}
\end{corollary}

{\bf Proof.}  It turns  out that in each of the first three cases $\Tor_j^R(R/\fa,\mathfrak{F}^i_{\fa}(M))$ is Artinian for all
$i$ and $j$, see respectively Theorems 3.2, 3.6 and  3.8.
Hence, Lemma 3.1 yields that in each of these cases $\Tor_j^R(R/\fm,\mathfrak{F}^i_{\fa}(M))$
is also Artinian for all $i$ and $j$, and so $\beta^j(\fm,\mathfrak{F}^i_{\fa}(M))
(:=\dim_{R/\fm}\Tor_j^R(R/\fm,\mathfrak{F}^i_{\fa}(M)))$ is
finite for all $i$ and $j$.

Now, assume that $R$ is regular of positive characteristic. We may
assume that $R$ is complete. Then by Lemma 2.1 ii), one has
$$\mathfrak{F}^i_{\fa}(M)\cong \Hom_R(H^{\dim R-i}_{\fa}(M,R),
E_R(R/\fm))$$ for all $i\geq 0$. This implies that $$\Tor_j^R(R/\fm,\mathfrak{F}^i_{\fa}(M))\cong  \Hom_R(\Ext_R^j(R/\fm,H^{\dim R-i}_{\fa}(M,R)),E_R(R/\fm))$$ for all $i$ and $j$. Thus
$\beta^j(\fm,\mathfrak{F}^i_{\fa}(M))$ is equal to the $j$th Bass number of
$H^{\dim R-i}_{\fa}(M,R)$, and so in the case iv), the claim follows by
\cite[Theorem 2.10]{DS}.
 $\Box$

\begin{example} i) The assumption $\dim R/\fa \leq 1$ cannot be dropped
in Theorem 3.6. To this end, let $k$ be a field, $R:=k[[W,X,Y,Z]]$ and
$\fa:=(W,X)\cap (Y,Z)$. Then $\dim R/\fa=2$. On the other hand, one has
$\mathfrak{F}^1_{\fa}(R)\cong R$, see \cite[Example 5.2]{Sch}. Hence
$\mathfrak{F}^1_{\fa}(R)/\fa \mathfrak{F}^1_{\fa}(R)$ is not Artinian.\\
ii) The assumption $\dim R\leq 2$ is really needed in Theorem 3.8 ii).
To realize this, let $k$ be a field, $R:=k[[X,Y,Z]]$ and
$\fa:=(XZ,YZ)$. Then by using Lemma 2.1 ii), we get the isomorphism
$\mathfrak{F}^1_{\fa}(R)/\fa \mathfrak{F}^1_{\fa}(R)\cong
\Hom_R(\Hom_R(R/\fa,H_{\fa}^2(R)),E_R(R/\fm))$. Since in view of the
proof of \cite[Theorem 2.2]{M1}, the module $\Hom_R(R/\fa,H_{\fa}^2(R))$
is not finitely generated, it follows that $\mathfrak{F}^1_{\fa}(R)/\fa \mathfrak{F}^1_{\fa}(R)$ is not Artinian.
\end{example}

\section{Formal grade and Cohen-Macaulayness}

Recall that an ideal $\fb$ of a ring $R$ is said to be {\it cohomologically
complete intersection} if $\cd_{\fb}(R)=\Ht \fb$, see \cite{HS}. Typical
examples of such ideals are set-theoretic complete intersection ideals.
Let $\fa$ be an ideal of a complete local ring $(R,\fm)$ and $M$ a finitely
generated $R$-module. Let $D_R^{\bullet}$ be a normalized dualizing complex
of $R$. Schenzel \cite[Theorem 3.5]{Sch} has proved the duality isomorphism
$\mathfrak{F}^{i}_{\fa}(M)\cong \Hom_R(H^{-i}_{\fa}(\Hom_R(M,D^{\bullet}_R)),E_R(R/\fm))$
for all $i$. This, in particular yields that $\fgrade(\fa,M)=-\sup\{i\in \mathbb{N}_0: H^i_{\fa}(\Hom_R(M,D^{\bullet}_R))\neq 0\}$.
In our first result in this section, we specialize this result in two particular situations.

\begin{theorem} Let $\fa$ be an ideal of a complete local ring
$(R,\fm)$ and $M$ a finitely generated $R$-module.
\begin{enumerate}
\item[i)] If $M$ is Cohen-Macaulay, then $\mathfrak{F}^i_{\fa}(M)\cong
\Hom_R(H_{\fa}^{\dim M-i}(K_M),E_R(R/\fm))$ for all $i\geq 0$. In particular,
$\fgrade(\fa,M)=\dim M-\cd_{\fa}(M)$.
\item[ii)] If $R$ is Cohen-Macaulay and $\fa$
cohomologically complete intersection, then $\mathfrak{F}^i_{\fa}(M)\cong
\Hom_R(\Ext^{\ell-i}_R(M,H^{\Ht \fa}_{\fa}(K_R)),E_R(R/\fm))$
for all $i\geq 0$, where $\ell:=\dim R/\fa$. In particular, $\fgrade(\fa,M)=
\dim R/\fa-\sup\{i:\Ext_R^i(M,H_{\fa}^{\Ht \fa}(K_{R}))
\neq 0\}$.
\end{enumerate}
\end{theorem}

{\bf Proof.}  Since $R$ is complete, it possesses a dualizing complex.
Let $D^{\bullet}_R$ be a normalized dualizing complex of $R$. Then by
\cite[Theorem 3.5]{Sch}, one has $$\mathfrak{F}^{i}_{\fa}(M)\cong
\Hom_R(H^{-i}_{\fa}(\Hom_R(M,D^{\bullet}_R)),E_R(R/\fm))$$ for all $i$.
By \cite [Proposition 2.4 b)]{Sch}, we have $$H_{\fm}^i(M)\cong \Hom_R(H^{-i}
(\Hom_R(M,D^{\bullet}_R)), E_R(R/\fm))$$ for all $i\geq 0$.
Assume that $M$ is Cohen-Macaulay. Then the two complexes
$\Hom_R(M,D^{\bullet}_R)$ and $K_M[\dim M]$ are quasi-isomorphic.
Hence $H^{-i}_{\fa}(\Hom_R(M,D^{\bullet}_R))\cong H_{\fa}^{\dim_RM-i}(K_M)$,
and so the first assertion of i) follows.

Now, assume that $R$ is Cohen-Macaulay and $\fa$ cohomologically complete
intersection. Since $M$ is finitely generated, one has ${\bf R}\Gamma_{\fa}
(\Hom_R(M,D^{\bullet}_R))\simeq {\bf R}\Hom_R(M,{\bf R}\Gamma_{\fa}(D^{\bullet}_R))$.
Using some of the standard properties of canonical modules, one can easily
deduce that $\grade(\fa,K_R)=\grade(\fa,R)$. On the other hand,
as $\Supp_RK_R=\Spec R$, \cite [Theorem 2.2]{DNT}
yields that $\cd_{\fa}(K_R)=\cd_{\fa}(R)$. So, because $\fa$ is cohomologically
complete intersection, it turns out that $H_{\fa}^i(K_R)=0$ for all $i\neq \Ht \fa$.
Hence, ${\bf R}\Gamma_{\fa}(K_R)\simeq H_{\fa}^{\Ht \fa}(K_R)[-\Ht \fa]$.
On the other
hand, by preceding paragraph, one has $D^{\bullet}_R\simeq K_R[\dim R]$. Thus
${\bf R}\Gamma_{\fa}(D^{\bullet}_R)\simeq
H^{\Ht \fa}_{\fa}(K_R)[\ell]$. Therefore, $${\bf R}\Gamma_{\fa}(\Hom_R(M,D^{\bullet}_R))
\simeq {\bf R}\Hom_R(M,{\bf R}\Gamma_{\fa}(D^{\bullet}_R))
\simeq {\bf R}\Hom_R(M,H^{\Ht \fa}_{\fa}(K_R)[\ell]),$$
and so
$$\begin{array}{ll}\mathfrak{F}^{i}_{\fa}(M)&\cong \Hom_R(H^{-i}({\bf R}\Gamma_{\fa}(\Hom_R(M,D^{\bullet}_R))),E_R(R/\fm))\\&\cong
\Hom_R(H^{-i}({\bf R}\Hom_R(M,H^{\Ht \fa}_{\fa}(K_R)[\ell]))),E_R(R/\fm))\\&\cong
\Hom_R(\Ext^{\ell-i}_R(M,H^{\Ht \fa}_{\fa}(K_R)),E_R(R/\fm))
\end{array}$$
for all $i$. This completes the proof of the first assertion of ii).

The proof of the last assertions of i) and ii) are similar to the proof of
the last assertion of Lemma 2.1 ii), and so we leave it to the reader.  $\Box$

The corollary below provides a new characterization of Cohen-Macaulay modules.

\begin{corollary} Let $(R,\fm)$ be a local ring and $M$ a nonzero finitely generated
$R$-module. Then the following are equivalent:
\begin{enumerate}
\item[i)] $M$ is Cohen-Macaulay.
\item[ii)] $\fgrade(\fa,M)+\cd_{\fa}(M)=\dim M$ for all ideals $\fa$ of $R$.
\item[iii)] $\fgrade(\fa,M)+\cd_{\fa}(M)=\depth M$ for all ideals $\fa$ of $R$.
\end{enumerate}
\end{corollary}

{\bf Proof.} By \cite[Proposition 3.3]{Sch}, we have $\mathfrak{F}^i_{\fa}(M)
\cong \mathfrak{F}^i_{\fa \widehat{R}}(\widehat{M})$ for all $i$. Hence,
without loss of generality, we may and do assume that $R$ is complete.

$i)\Rightarrow ii)$ This is immediate by Theorem 4.1 i).

$ii)\Rightarrow iii)$ Consider the ideal $\fa:=0$. Then
$\mathfrak{F}_{\fa}^i(M)\cong H^{i}_{\fm}(M)$ for all $i$.
So,  $\fgrade(\fa,M)=\depth_RM$ and $\cd_{\fa}(M)=0$. Thus
$\dim M=\depth_RM$. This yields iii).

$iii)\Rightarrow i)$  Let $\fa:=\fm$. Then
$\mathfrak{F}_{\fa}^0(M)\cong M$, and so $\fgrade(\fa,M)=0$. On the
other hand, Grothendieck's
non-vanishing Theorem asserts that $\cd_{\fa}(M)=\dim M$. Thus
$\dim M=\depth_RM$, as required. $\Box$

Let $M$ be a Cohen-Macaulay module over a local ring $(R,\fm)$. By Corollary 4.2,
we know that  $\fgrade(\fa,M)+\cd_{\fa}(M)=\dim M$ for all ideals $\fa$ of $R$.
It would be interesting to know whether the same equality remains true for some
special types of ideals and of modules. In view of Theorems 3.2
and 3.6, principal and one dimensional ideals might be appropriate
candidates for our desired ideals. Also,  some generalizations of
the notion of Cohen-Macaulay modules could be appropriate candidates
for our desired modules. The following examples indicate that for these
types of ideals the above equality does not hold even for sequentially
Cohen-Macaulay modules, Buchsbaum rings, quasi-Gorenstein rings
and approximately Cohen-Macaulay rings.

\begin{example} i) Let $(R,\fm)$ be a 2-dimensional regular local ring and
$\fa$ an ideal of $R$ with $\dim R/\fa=1$. The Hartshorne-Lichtenbaum
Vanishing Theorem yields that $\cd_{\fa}(R)=1$ and clearly  $\cd_{\fa}(R/\fm)=0$.
Hence, Corollary 4.2 implies that $\fgrade(\fa,R)=1$ and $\fgrade(\fa,R/\fm)=0$.
Set $M:=R\oplus R/\fm$. Then $M$ is a 2-dimensional sequentially Cohen-Macaulay
$R$-module. We have $\cd_{\fa}(M)=1$ and $$\fgrade(\fa,M)=\min\{\fgrade(\fa,R),
\fgrade(\fa,R/\fm)\}=0.$$ Hence $\fgrade(\fa,M)+\cd_{\fa}(M)<\dim M$.

ii) Let $X$ be an Abelian variety and $R:=\bigoplus_{n\in\mathbb{Z}}
H^0(X,L^{\otimes n})$, where $L$ is a very ample invertible sheaf on $X$. Let $\fm:=
\bigoplus_{n>0} H^0(X, L^{\otimes n})$ and assume that $g:=\dim X>0$.
Then $\depth R_{\fm}=2$, $\dim R_{\fm}=g+1$,  $H^{i}_{\fm}(R_{\fm})$ is a
finitely generated nonzero $R_{\fm}$-module for all $2\leq i\leq g$ and
$H^{g+1}_{\fm}(R_{\fm})\cong E_{R_{\fm}}(R_{\fm}/\fm R_{\fm})$. (This
example is due to Schenzel, see \cite[page 235]{SV} for more details.)
So, $R_{\fm}$ is a Buchsbaum quasi-Gorenstein local ring.
Now, take $g\geq 3$ and let $\fa$ be a nonzero principal ideal of $R_{\fm}$.
Then, by \cite[Theorem 4.9]{Sch}, one has
$$\fgrade(\fa,R_{\fm})=\inf\{i-\cd_{\fa}(K^i_{R_{\fm}}):i=0,\cdots, g+1\}=2.$$
Hence, since $\cd_{\fa}(R_{\fm})\leq 1$, we deduce that $\fgrade(\fa,R_{\fm})+
\cd_{\fa}(R_{\fm})<\dim R_{\fm}$.

iii) Let $k$ be a field. Consider the $2$-dimensional complete local ring
$R:=\frac{k[[X,Y,Z]]}{(X)\cap(Y,Z)}$. One can check that $R$ is approximately
Cohen-Macaulay (i.e. there exists an element $a$ of $R$ such that
$R/a^nR$ is a Cohen-Macaulay ring of dimension $1$ for every integer $n>0$).
Set $\fa:=(x)$. Note that by \cite [Theorem 4.12]{Sch}, we have
$\fgrade(\fa,R)\leq \dim(R/\fa+\fp)$
for all $\fp\in \Ass_RR$. Now, by applying this result to $\fp:=(y,z)$, we
find that $\fgrade(\fa,R)=0$. On the other hand, $\cd_{\fa}(R)\leq 1$, and so
$\fgrade(\fa,R)+\cd_{\fa}(R)<\dim R.$

iv) One can restate the equivalence $i)\Leftrightarrow ii)$ in Corollary 4.2
by saying that $R$ is Cohen-Macaulay if and only if for every nilpotent
ideal $\fa$ of $R$, $\fgrade(\fa,R)+\cd_{\fa}(R)=\dim R$. Now, we give an
example of a local ring $(R,\fm)$ such that for any ideal $\fa$ of $R$, the
formula $\fgrade(\fa,R)+\cd_{\fa}(R)=\dim R$ holds if and only if $\fa$ is
non-nilpotent. To this end, let $k$ be a field and $R:=k[[X,Y]]/(XY,Y^2)$.
Let $\fa$ be a nilpotent ideal of $R$. Then $\mathfrak{F}^0_{\fa}(R)=
\underset{n}{\vpl}H^0_{\fm}(R/\fa^n
 R)=H^0_{\fm}(R)\neq 0$. Hence, $\fgrade(\fa,R)=\cd_{\fa}(R)=0$, and so
$$\fgrade(\fa,R)+\cd_{\fa}(R)<1=\dim R.$$ Next, let $\fa$ be a non-nilpotent
ideal of $R$. Then $\fa$ is an $\fm$-primary ideal of $R$. Hence,
$\cd_{\fa}(R)=1$ and $\fgrade(\fa,R)=0$, and so $\fgrade(\fa,R)+
\cd_{\fa}(R)=\dim R$.
\end{example}

Let $M$ be an $R$-module and $L$ a pure submodule of $M$. Let $\fa$ be
an ideal of $R$ such that the equality $\fgrade(\fa,M)+\cd_{\fa}(M)=\dim M$
holds. In the next two propositions, we investigate the question whether this
equality implies that $\fgrade(\fa,L)+\cd_{\fa}(L)=\dim L$.

\begin{proposition} Let $\fa$ be an ideal of a local ring $(R,\fm)$ and $M$ a
finitely generated $R$-module. Assume that $L$ is a pure submodule of $M$.
\begin{enumerate}
\item[i)] $\fgrade(\fa,L)\geq \fgrade(\fa,M)$.  In particular, if $M$
is Cohen-Macaulay, then $L$ is also Cohen-Macaulay.
\item[ii)] If $M$ and $L$ have the same support, then $\fgrade(\fa,M)+
\cd_{\fa}(M)=\dim M$ implies that $\fgrade(\fa,L)+\cd_{\fa}(L)=\dim L$,
and so $\fgrade(\fa,L)=\fgrade(\fa,M)$.
\end{enumerate}
\end{proposition}

{\bf Proof.} i) Since $L$ is a pure submodule of $M$, one concludes that the
natural map $L/\fa^n L\lo M/\fa^n M$ is pure for all $n\geq 0$. Now, \cite [Corollary 3.2 a)]{Ke}
implies that the induced map $$H_{\fm}^i(L/\fa^n L)\lo H_{\fm}^i(M/\fa^n M)$$ is
injective for all $i$ and $n$. For each $i$,  the inverse system
$\{H_{\fm}^i(L/\fa^n L)\}_{n\in \mathbb{N}}$ satisfies the Mittag-Leffler
condition. Thus,
it follows that $\mathfrak{F}_{\fa}^i(L)$ is isomorphic to a submodule of
$\mathfrak{F}_{\fa}^i(M)$, and so $\fgrade(\fa,L)\geq \fgrade(\fa,M)$. For the
remaining assertion of i), note that for any finitely generated $R$-module $N$,
one has $\fgrade(0,N)=\depth N$. Hence by applying the first assertion of i), for
the zero ideal, we have $$\depth L=\fgrade(0,L)\geq \fgrade(0,M)=\depth M.$$ So, the equality $\depth M=\dim M$ yields the equality $\depth L=\dim L$.

ii) Assume that $\Supp_RL=\Supp_RM$. Then $\dim L=\dim M$ and \cite [Theorem 2.2]{DNT}
implies that $\cd_{\fa}(L)=\cd_{\fa}(M)$. Suppose that $\fgrade(\fa,M)+\cd_{\fa}(M)=\dim M$.
Then by i) and \cite [Corollary 4.11]{Sch}, we have
$$\dim L\geq \fgrade(\fa,L)+\cd_{\fa}(L)\geq \fgrade(\fa,M)+\cd_{\fa}(M)=\dim M.$$
This finishes the proof of ii). $\Box$

Let $G$ be a group of automorphisms of $R$. Recall that the ring of invariants $R^G$
is defined to be the set of all elements of $R$, which are invariant under the
action of $G$. Our next result can be considered as a slight generalization of the Hochster-Eagon result on Cohen-Macaulayness of invariant rings which it corresponds
to the case $\fb=0$.

\begin{proposition} Let $(R,\fm)$ be a local ring and $G$ a group of
automorphisms of $R$ such that $R$ is integral over $R^{G}$. Assume that
there exists a Reynolds operator $\rho:R\lo R^{G}$. Let
$\fb$ be an ideal of $R^G$ such that $\fgrade(\fb R,R)+\cd_{\fb R}(R)=\dim R$. Then
$\fgrade(\fb,R^G)=\fgrade(\fb R,R)$ and $\cd_{\fb}(R^G)=\cd_{\fb R}(R)$,
in particular, $\fgrade(\fb,R^G)+\cd_{\fb}(R^G)=\dim R^G$.
\end{proposition}

{\bf Proof.} Since $\rho|_{R^G}=id_{R^G}$, one has $R=R^G\bigoplus X$ for
some $R^G$-module $X$. It follows easily that $R^G$ is also a Noetherian ring.
Because $R$ is integral over $R^G$, it turns out that $R^G$ is also local and
$\dim R=\dim R^G$. By the Independence Theorem for local cohomology modules,
we have $$H^i_{\fb R}(R)\cong H^i_{\fb}(R)\cong H^i_{\fb}(R^G)\oplus
H^i_{\fb}(X),$$ and so $\cd_{\fb}(R^G)\leq \cd_{\fb}(R)=\cd_{\fb R}(R)$.
Now, since $\cd_{\fb}(R^G)$ is the supremum of $\cd_{\fb}(M)$'s, where
$M$ runs over all $R^G$-modules, one concludes that $\cd_{\fb}(R^G)=
\cd_{\fb R}(R)$. On the other hand, by using Lemma 2.1 i), one has
$$\mathfrak{F}^i_{\fb R}(R)\cong
\mathfrak{F}^i_{\fb}(R)\cong\mathfrak{F}^i_{\fb}(R^G)\oplus
\mathfrak{F}^i_{\fb}(X),$$ and so $\fgrade(\fb R,R)\leq \fgrade(\fb,R^G)$.
Therefore,
$$\fgrade(\fb,R^G)+\cd_{\fb}(R^G)\geq\fgrade(\fb R,R)+\cd_{\fb R}(R)=\dim
R=\dim R^G.$$ The reverse of the above inequality always holds by
\cite[Corollary 4.11]{Sch}. So, $$\fgrade(\fb,R^G)=\fgrade(\fb R,R).$$
This completes the proof. $\Box$

Schenzel \cite[Corollary 4.11]{Sch} has proved that $\fgrade(\fa,M)\leq
\dim M-\cd_{\fa}(M)$. In the next result, we establish a lower bound for
$\fgrade(\fa,M)$.

\begin{theorem} Let $\fa$ be an ideal of a local ring $(R,\fm)$ and $M$ a
finitely generated $R$-module. Then
$$\depth M -\cd_{\fa}(M)\leq \fgrade(\fa,M)\leq \dim M-\cd_{\fa}(M).$$
\end{theorem}

{\bf Proof.} The right-hand inequality holds by \cite[Corollary 4.11]{Sch}.
By \cite[Proposition 3.3]{Sch}, we have $\mathfrak{F}^i_{\fa}(M)
\cong \mathfrak{F}^i_{\fa \widehat{R}}(\widehat{M})$ for all $i$. Therefore,
without loss of generality, we may and do assume that $R$ is complete.
Then, by Cohen's Structure Theorem  $R$ is a homomorphic image of a regular
complete local ring $(T,\fn)$ 'say. So, $R\cong T/J$ for some ideal $J$ of
$T$. Set $b:=\fa\cap T$. Then by Lemma 2.1, it follows that
$$\mathfrak{F}^i_{\fa}(M)\cong \mathfrak{F}^i_{\fb}(M)\cong
\Hom_T(H^{\dim T-i}_{\fb}(M,T),E_T(T/\fn))$$ for all $i$, and $\fgrade(\fa,M)=\dim T-
\cd_{\fb}(M,T)$. By using \cite[Corollary 2.10]{DH}, one has
$$\cd_{\fb}(M,T)\leq \pd_TM+\cd_{\fb}(M\otimes_TT).$$ Hence, by
the Auslander-Buchsbaum formula and the Independence Theorem for local
cohomology modules, we deduce that
$$
\begin{array}{ll} \fgrade(\fa,M)&\geq \dim T-\pd_TM-\cd_{\fb}(M\otimes_TT)
\\&=\depth_TM-\cd_{\fb}(M)\\&=\depth_RM-\cd_{\fa}(M).  \ \  \Box\\
\end{array}
$$

\begin{remark} i) Let $\fa$ be an ideal of a local ring $(R,\fm)$ and $L$
and $M$ two finitely generated $R$-modules such that $\Supp_RL\subseteq \Supp_RM$.
Then by \cite [Theorem 2.2]{DNT}, we know that $\cd_{\fa}(L)\leq \cd_{\fa}(M)$.
In particular, one has $\cd_{\fa}(L)=\cd_{\fa}(M)$, whenever $\Supp_RL=\Supp_RM$.
One might expect that the assumption $\Supp_RL=\Supp_RM$ forces the equality
$\fgrade(\fa,L)=\fgrade(\fa,M)$. But, as it is clear by
\cite [Example 4.10]{Sch} (or Example 4.3 i)),
this is not the case in general. However, if $L$ is a pure submodule of $M$
such that $\Supp_RL=\Supp_RM$ and  $\fgrade(\fa,M)+\cd_{\fa}(M)=\dim M$,
then Proposition 4.4 ii) implies that $\fgrade(\fa,L)=\fgrade(\fa,M)$.
Note that the assumption $\fgrade(\fa,M)+\cd_{\fa}(M)=\dim M$ is really
needed. To see this, let $M$ and $L:=R$ be as in Example 4.3 i).

ii) Let $(R,\fm)$ be a local ring and $G$ a group of automorphisms of $R$.
For each $r\in R$, the orbit of $r$ under the action of $G$ is denoted by
$G_r$. The group $G$ is said to be {\it locally finite} if for each element
$r\in R$, the set $G_r$ is finite. Let $G$ be a locally finite group of
automorphisms of $R$ such that $|G_r|$ is a unit in $R$ for every $r\in R$
(e.g. $G$ is a finite group such that $|G|$ is a unit in $R$).
Then the map $\rho:R \lo R^G$ given by the assignment $r$ to
$\frac{1}{|G_r|}\underset{s\in G_r}\Sigma s$ is a Reynolds operator and
$R$ is integral over $R^{G}$. So by Proposition 4.5, if for an ideal $\fb$ of
$R^G$, one has $\fgrade(\fb R,R)+\cd_{\fb R}(R)=\dim R$, then
$\fgrade(\fb,R^G)+\cd_{\fb}(R^G)=\dim R^G$.

iii) Note that Corollary 4.2 can also be deduced easily from Theorem 4.6.
\end{remark}

\begin{acknowledgement}  In an earlier version, we have used a spectral
sequence argument for the proof of Theorem 4.1. Professor Peter Schenzel
has pointed out to us that they can also be deduced from \cite [Theorem 3.5]{Sch}.
Here, we exposed his proposed
argument in place of our original one. The authors would like to express
their thanks to him for his valuable comments.
\end{acknowledgement}



\begin{thebibliography}{99}

\bibitem[ADT]{ADT}{ M. Asgharzadeh, K. Divaani-Aazar and M. Tousi},
{\it Finiteness dimension of local cohomology modules and its dual notion},
J. Pure Appl. Algebra, {\bf 213}(3), (2009), 321-328.

\bibitem[B]{B}{M.H. Bijan-Zadeh}, {\it A common generalization
of local cohomology theories}, Glasgow Math. J., {\bf 21}(2),
(1980), 173-181.

\bibitem[BS]{BS}{M. Brodmann and R.Y. Sharp}, {\it Local cohomology:
an algebraic introduction with geometric applications}, Cambridge Univ.
Press, {\bf 60}, Cambridge, (1998).

\bibitem[DM]{DM}{D. Delfino and T. Marley}, {\it Cofinite modules and local
cohomology}, J. Pure Appl. Algebra, {\bf 121}(1), (1997), 45-52.

\bibitem[DH]{DH}{K. Divaani-Aazar and  A. Hajikarimi}, {\it Generalized local
cohomology modules and homological Gorenstein dimensions}, Comm. Algebra, to appear.

\bibitem[DS]{DS}{K. Divaani-Aazar and R. Sazeedeh}, {\it Cofiniteness of
generalized local cohomology modules}, Colloq. Math., {\bf 99}(2), (2004),
283-290.

\bibitem[DNT]{DNT}{K. Divaani-Aazar and R. Naghipour and M. Tousi}, {\it
Cohomological dimension of certain algebraic varieties}, Proc. Amer. Math. Soc.,
{\bf 130}(12), (2002), 3537-3544.

\bibitem[DT]{DT}{K. Divaani-Aazar and M. Tousi}, {\it Some remarks on
coassociated primes}, J. Korean Math. Soc., {\bf 36}(5), (1999), 847-853.

\bibitem[Hel1]{Hel1}{M. Hellus},  {\it On the associated primes of Matlis duals
of top local cohomology modules}, Comm. Algebra, {\bf 33}(11),
(2005), 3997-4009.

\bibitem[Hel2]{Hel2}{M. Hellus}, {\it A note on the injective dimension of
local cohomology modules},  Proc. Amer. Math. Soc., {\bf 136}(7), (2008), 2313-2321.

\bibitem[Hel3]{Hel3}{M. Hellus}, {\it Local Cohomology and
Matlis duality}, Habilitationsschrift, Universit\"{a}t Leipzig, (2007),
arXiv:0703124.

\bibitem[HS]{HS}{M. Hellus and P. Schenzel},  {\it On cohomologically complete
intersections}, J. Algebra, {\bf 320}(10), (2008), 3733-3748.

\bibitem[Her]{Her}{J. Herzog}, {\it Komplexe Aufl\"{o}sungen und Dualit\"{a}t in der
lokalen Algebra}, Habilitationsschrift, Universit\"{a}t Regensburg, (1974).

\bibitem[HE]{HE}{M. Hochster and J.A. Eagon}, {\it Cohen-Macaulay rings, invariant
theory, and the generic perfection of determinantal loci}, Amer. J. Math., {\bf 93},
(1971), 1020-1058.

\bibitem[HK]{HK}{C. Huneke and J. Koh}, {\it Cofiniteness and vanishing of local
cohomology modules}, Math. Proc. Cambridge Philos. Soc., {\bf 110},
(1991), 421-429.

\bibitem[HS]{HS}{C. Huneke and R.Y. Sharp}, {\it Bass numbers of local cohomology
modules}, Trans. Amer. Math. Soc., {\bf 339}(2), (1993), 765-779.

\bibitem[I]{I}{L. Illusie}, {\it Grothendieck's existence theorem in formal
geometry}, in Fundamental algebraic geometry, Mathematical Surveys and Monographs,
{\bf 123}, American Mathematical Society, Providence, RI, (2005), 179-234.

\bibitem[K]{K}{K.I. Kawasaki}, {\it Cofiniteness of local cohomology modules for
principal ideals}, Bull. London Math. Soc., {\bf 30}(3), (1998), 241-246.

\bibitem[Ke]{Ke}{G. Kempf}, {\it The Hochster-Roberts theorem of invariant
theory}, Michigan Math. J., {\bf 26}(1), (1979),  19-32.

\bibitem[L]{L}{G. Lyubeznik}, {\it Finiteness properties of local cohomology
modules (an application of $D$-modules to commutative algebra)}, Invent.
Math., {\bf 113}(1), (1993), 41-55.

\bibitem[M1]{M1}{L. Melkersson}, {\it Properties of cofinite modules and applications
to local cohomology}, Math. Proc. Cambridge Philos. Soc., {\bf 125}(3), (1999), 417-423.

\bibitem[M2]{M2}{L. Melkersson}, {\it Cohomological properties of modules with
secondary representations}, Math. Scand., {\bf 77}(2), (1995), 197-208.

\bibitem[O]{O}{A. Ogus}, {\it Local cohomological dimension of algebraic varieties},
Ann. Math., {\bf 98}, (1973), 327-365.

\bibitem[PS]{PS}{C. Peskine and L. Szpiro}, {\it Dimension projective finie
et cohomologie locale}, Publ. Math. I.H.E.S., {\bf42}, (1972), 47-119.

\bibitem[R]{R}{J. Rotman}, {\it An Introduction to Homological
Algebra}, Academic Press, San Diego, (1979).

\bibitem[Sch]{Sch}{P. Schenzel},  {\it On formal  local cohomology and
connectedness}, J. Algebra, {\bf 315}(2), (2007), 894-923.

\bibitem[SV]{SV}{J. St$\ddot{u}$ckrad and W. Vogel}, {\it Buchsbaum rings
and applications}, Springer-Verlag, Berlin, (1986).

\bibitem[V]{V}{W.V. Vasconcelos}, {\it Divisor theory in module categories},
North-Holland Mathematics Studies, {\bf 14}, (1974).

\bibitem[Y]{Y}{K.I. Yoshida}, {\it Cofiniteness of local cohomology modules for
ideals of dimension one}, Nagoya Math. J., {\bf 147}, (1997), 179-191.

\bibitem[Z]{Z}{H. Z\"{o}schinger}, {\it Der Krullsche Durchschnittssatz f\"{u}r
kleine Untermoduln}, Arch. Math. (Basel), {\bf 62}(4), (1994), 292-299.

\end{thebibliography}
\end{document}